\documentclass[11pt]{article}
\usepackage{amsfonts,amsmath,amssymb,amsthm}
\usepackage{graphicx,url}
\usepackage{natbib}
\usepackage{flexisym}
\usepackage[lined,boxed,linesnumbered]{algorithm2e}

\newtheorem{proposition}{Proposition}[section]
\newtheorem{example}{Example}[section]
\newtheorem{definition}{Definition}[section]
\newtheorem{lemma}{Lemma}[section]
\newtheorem{remark}{Remark}[section]

\def\sgn{\mathop{\rm sgn}\nolimits}

\def\tr{\mathop{\rm tr}\nolimits}

\def\dist{\mathop{\rm dist}\nolimits}

\def\dom{\mathop{\rm dom}\nolimits}
\def\argmin{\mathop{\rm argmin}\nolimits}
\def\amp{\mathop{\;\:}\nolimits}
\newcommand{\ba}{\boldsymbol{a}}
\newcommand{\bb}{\boldsymbol{b}}

\newcommand{\bs}{\boldsymbol{s}}

\newcommand{\bu}{\boldsymbol{u}}
\newcommand{\bv}{\boldsymbol{v}}

\newcommand{\bx}{\boldsymbol{x}}
\newcommand{\by}{\boldsymbol{y}}
\newcommand{\bz}{\boldsymbol{z}}
\newcommand{\bA}{\boldsymbol{A}}
\newcommand{\bB}{\boldsymbol{B}}

\newcommand{\bD}{\boldsymbol{D}}

\newcommand{\bI}{\boldsymbol{I}}
\newcommand{\bJ}{\boldsymbol{J}}
\newcommand{\bK}{\boldsymbol{K}}

\newcommand{\bM}{\boldsymbol{M}}

\newcommand{\bO}{\boldsymbol{O}}
\newcommand{\bP}{\boldsymbol{P}}
\newcommand{\bQ}{\boldsymbol{Q}}

\title{Nonconvex Optimization via MM Algorithms: Convergence Theory}

\author{Kenneth Lange\\
Departments of Computational Medicine, \\
Human Genetics, and Statistics\\
University of California at Los Angeles\\
\url{klange@ucla.edu} \\ \\
Joong-Ho Won\\
Department of Statistics\\
Seoul National University\\
\url{wonj@stats.snu.ac.kr} \\ \\
Alfonso Landeros\\
Department of Computational Medicine\\
University of California at Los Angeles\\
\url{alanderos@ucla.edu} \\ \\
Hua Zhou\\
Department of Biostatistics\\
University of California at Los Angeles\\
\url{huazhou@ucla.edu}}

\date{\today}

\begin{document}
\maketitle

\newpage

\begin{abstract}
The majorization-minimization (MM) principle is an extremely general framework for deriving optimization algorithms. It includes the expectation-maximization (EM) algorithm, proximal gradient algor\-ithm, concave-convex procedure, quadratic lower bound algorithm, and proximal distance algorithm  as special cases. Besides numerous applications in statistics, optimization, and imaging, the MM principle finds wide applications in large scale machine learning problems such as matrix completion, discriminant analysis, and nonnegative matrix factorizations. When applied to nonconvex optimization problems, MM algorithms enjoy the advantages of convexifying the objective function, separating variables, numerical stability, and ease of implementation. However, compared to the large body of literature on other optimization algorithms, the convergence analysis of MM algorithms is  scattered and problem specific. This survey presents a unified treatment of the convergence of MM algorithms. With modern applications in mind, the results encompass non-smooth objective functions and non-asymptotic analysis.
\end{abstract}

\section{Background}

The majorization-minimization (MM) principle for constructing optimization algorithms \citep{BeckerYangLange97MM,Lange00OptTrans,hunter2004tutorial} finds broad range of applications in
\begin{itemize}
\item statistics: multidimensional scaling \citep{BorgGroenen2005modern}, quantile regression \citep{HunterLange00QuantileRegMM}, ranking sports teams \citep{Hunter04BreadleyTerry}, variable selection \citep{HunterLi05MM,Yen11MMSpikeSlab,BienTibshirani11SparseCov,LeeHuang13MMSpPCA}, multivariate distributions \citep{ZhouLange10DMMLE,ZhangZhouZhouSun17mglm}, variance components models \citep{ZhouHuZhouLange19VCMM}, robust covariance estimation \citep{SunBabuPalomar15RobustCov}, and survival models \citep{HunterLange02PropOdds,DingTianYuen15CoxMM};

\item optimization: geometric and sigmoid programming \citep{LangeZhou14GPMM},  proximal distance algorithm \citep{chi2014distance,XuChiLange17GLMProxDist,KeysZhouLange19ProxDist};

\item imaging: transmission and positron tomography \citep{Lange84PET}, wavelets \citep{Figueiredo07MMWavelet}, magnetic resonance imaging, sparse deconvolution;

\item and machine learning: nonnegative matrix factorization \citep{LeeSeung99NNMF}, matrix completion \citep{Mazumder10SVDReg,chi2013genotype}, clustering \citep{ChiLange13ConvexClustering,XuLange19Kmeans}, discriminant analysis \citep{WuLange10DA}, support vector machines \citep{Nguyen17MMIntro}.
\end{itemize}
The recent book \citep{Lange16MMBook} and survey papers \citep{SunBabuPalomar17MM,Nguyen17MMReview} give a comprehensive overview of MM algorithms. 

The MM principle involves majorizing the objective function $f(\bx)$ by a surrogate function $g(\bx \mid \bx_n)$ around the current iterate $\bx_n$ of a search. Majorization is defined by the two conditions
\begin{eqnarray}
	f(\bx_n) &=& g(\bx_n \mid \bx_n) \label{eqn:tangency}\\
	f(\bx) &\le& g(\bx \mid \bx_n), \quad \bx \ne \bx_n. \label{eqn:dominance}
\end{eqnarray}
In other words, the surface $\bx \mapsto g(\bx \mid \bx_n)$ lies above the surface $\bx \mapsto f(\bx)$ and is tangent to it at the point $\bx = \bx_n$. Construction of the majorizing function $g(\bx \mid \bx_n)$ constitutes the first M of the MM algorithm. The second M of the algorithm minimizes the surrogate $g(\bx \mid \bx_n)$ rather than $f(\bx)$. If $\bx_{n+1}$ denotes the minimizer of $g(\bx \mid \bx_n)$, then this action forces the descent property $f(\bx_{n+1}) \le f(\bx_n)$. This fact follows from the inequalities
\begin{eqnarray*}
	f(\bx_{n+1}) \le g (\bx_{n+1} \mid \bx_n) \le g(\bx_n \mid \bx_n) = f(\bx_n),
\end{eqnarray*}
reflecting the definition of $\bx_{n+1}$ and the tangency condition. 

The same principle applied to the maximization problems leads to the  minorization-maximization algorithms that monotonically increase the objective values. The celebrated EM algorithm in statistics is a special case of the minorization-maximization algorithm as the E-step constructs a Q-function that satisfies the minorization properties. Derivation of EM algorithm hinges upon the notion of missing data and conditional expectation while that of MM algorithm upon clever use of inequalities. For most problems where a EM algorithm exists, the MM derivation often leads to the same algorithm. Notable exceptions include the maximum likelihood estimation (MLE) of the Dirichlet-Multinomial model \citep{ZhouLange10DMMLE,ZhouZhang12EMvsMM} and the variance components model \citep{ZhouHuZhouLange19VCMM}. However the MM principle has much wider applications as it applies to both minimization and maximization problems and does not rely on the notion of missing data. 

\section{Convergence Theorems}

Throughout, we denote by $\mathcal{X} \subset \mathbb{R}^d$ the subset underlying our problems. All of the functions we consider have domain $\mathcal{X}$ and are extended real-valued with range $\mathbb{R} \cup \{\infty\}$. The interior of set $S$ is denoted by $\mathbf{int}S$, and its closure by $\mathbf{cl}S$.

The following concepts are useful.
\begin{definition}[Effective domain]
	The effective domain of a function $f$ is defined and denoted by
	$$
	\dom{f} = \{\bx \in \mathcal{X}: f(\bx) < \infty \}.
	$$
\end{definition}
\begin{definition}[Properness]
	Function $f(\bx)$ is called \emph{proper} if $\dom{f} \neq \emptyset$.
\end{definition}
\begin{definition}[Directional derivatives]
    The directional derivative of function $f$ at $\bx\in\mathcal{X}$ is defined and denoted as
    $$
    d_{\bv}f(\bx) = \lim_{t\downarrow 0} \frac{f(\bx+t\bv)-f(\bx)}{t}
    $$
    if the limit exists.
\end{definition}
\noindent If $f$ is differentiable at $\bx$, then $d_{\bv}f(\bx)=\langle \nabla f(\bx), \bv \rangle$.
\begin{definition}[$L$-smoothness]
	Function $f$ is said to be $L$-smooth with respect to a norm $\|\cdot\|$ if it is differentiable on $\mathbf{int}\dom{f}$ and the gradient $\nabla f$ is Lipschitz continuous with a Lipschitz constant $L$:
	$$
	\|\nabla f(\bx) - \nabla f(\by)\| \le L\|\bx - \by\|,
	\quad
	\forall \bx, \by \in \mathbf{int}\dom{f}
	.
	$$
\end{definition}
\noindent It can be shown that $f(\bx)$ is $L$-smooth if and only if
$$
	f(\bx) \le f(\by) + \langle \nabla f(\by), \bx - \by \rangle + \frac{L}{2}\|\bx - \by\|^2 
	, \quad
	\forall \bx, \by \in \mathbf{int}\dom{f}.
$$
\begin{definition}[Strong convexity]
	Function $f$ is called $\mu$-strongly convex with respect to a norm $\|\cdot\|$, $\mu\ge 0$, if $f(\bx) - \frac{\mu}{2}\|\bx\|^2$ is convex.
\end{definition}
\noindent It can be shown that if $f(\bx)$ is $\mu$-strongly convex and has its minimum at $\by$, then
$$
	f(\bx) - f(\by) \ge \frac{\mu}{2}\|\bx - \by\|^2
	.
$$
\begin{definition}[Tangent vector, tangent cone]
    For a closed nonempty set $C\subset\mathcal{X}$, 
    the \emph{tangent cone} of $C$ at $\bx$ is 
    $$
    T_C(\bx) = \big\{\bv \in \mathcal{X} : \exists \{\bx_n\} \subset C , \{t_n\} \subset \mathbb{R}  \text{ such that } t_n \downarrow 0, \bx_n \to \bx \text{ and } \frac{\bx_n-\bx}{t_n} \to \bv\big\},
    $$
    where the notation $t_n\downarrow 0$ means that $t_n$ approaches 0 from above.
    A vector $\bv \in T_C(\bx)$ is said to be a \emph{tangent vector} of $C$ at $\bx$.
\end{definition}

\subsection{Classical convergence theorem}

Consider the problem of minimizing the objective function $f$ over a closed nonempty set $C\subset\mathcal{X}$. The following is immediate from the decent property of the MM algorithms:
\begin{proposition} 
    Let $\{\bx_n\}\subset\mathcal{X}$ be the iterates generated by an MM algorithm.
	Assume (a) $\bx_n \in C$ for each $n$.
	Then the sequence of objective values $\{f(\bx_n)\}$ monotonically decreases. 
	Furthermore, if (b) $p^\star = \inf_{x\in C}f(\bx) > -\infty$, then $\{f(\bx_n)\}$ converges.
\end{proposition}

Whether the limit is the desired minimum and whether the iterate $\{\bx_n\}$ will converge to a minimizer is more subtle.
For the latter,
a classical theory of convergence in nonlinear optimization algorithms is due to \cite{ZangwillMond69Book}. We first recap Zangwill's theory following the modern treatment of \cite{LuenbergerYe08Book}. Note that most of the iterative optimization algorithms, including the MM algorithms, generate a sequence $\{\bx_n\}$ by mapping $\bx_n\in\mathcal{X}$ to another point $\bx_{n+1}\in\mathcal{X}$. For example, in MM algorithms, $\bx_{n+1}$ is a point that minimizes the surrogate function $g(\bx|\bx_n)$ in $\mathcal{X}$. However, such a minimizer may not be unique unless the $g(\bx|\bx_n)$ satisfies certain assumptions. Rather, $\bx_{n+1}$ is \emph{one of the minimizers} of $g(\bx|\bx_n)$ and can be written as $\bx_{n+1}\in\argmin_{\bx\in C}g(\bx|\bx_n)$. Thus we may in general define an algorithm map as a \emph{set-valued map}:
\begin{definition}[Algorithm map]
    An \emph{algorithm map} $M$ is a mapping defined on $\mathcal{X}$ that assigns to every point $\bx\in\mathcal{X}$ a subset of $\mathcal{X}$.
\end{definition}
\noindent Among which point  $M(\bx_n)$ to choose as $\bx_{n+1}$ depends on the specific details of the actual optimization algorithm. If $M$ is a single-valued map, i.e., $M(\bx)$ is singleton for all $\bx\in\mathcal{X}$, we write $\bx_{n+1}=M(\bx_n)$.

A desirable property of an algorithm map is closure, which extends continuity of single-valued maps to set-valued ones:
\begin{definition}[Closure]
    A set-valued map $M$ from $\mathcal{X}$ to $\mathcal{Y}$ is said to be \emph{closed} at $\bx\in\mathcal{X}$ if $\by\in M(\bx)$ whenever $\{\bx_n\}\subset\mathcal{X}$ converges to $\bx$ and $\{\by_n: \by_n\in\bM(\bx_n)\}$ converges to $\by$. 
    The map $M$ is said to be closed on $\mathcal{X}$ if it is closed at each point of $\mathcal{X}$.
\end{definition}

The celebrated Zangwill's global convergence theorem is phrased in terms of an algorithm map $M$, a solution set $\Gamma$, and a \emph{descent function} $u$:
\begin{lemma}[Convergence Theorem A, \cite{ZangwillMond69Book}]\label{lem:zangwill}
Let the point-to-set map $M: \mathcal{X} \to \mathcal{X}$
determine an algorithm that given a point
$\bx_0 \in \mathcal{X}$ generates the sequence $\{\bx_n\}$.
Also let a solution set $\Gamma \subset \mathcal{X}$ be given.
Suppose that: 
    \begin{enumerate}
        \item all points $\bx_n$ are in a compact set $C \subset \mathcal{X}$;
        \item there is a continuous function $u : \mathcal{X} \to \mathbb{R}$ such that (a) if $\bx \notin \Gamma$, $u(\by) < u(\bx)$ for all $\by \in M(\bx)$, and (b) if $\bx \in \Gamma$, then either the algorithm terminates or  $u(\by) \le u(\bx)$ for all $\by \in M(\bx)$;
        \item the map $M$ is closed at $\bx$ if $\bx \notin \Gamma$.
    \end{enumerate}
Then either the algorithm stops at a solution, or the limit of any convergent subsequence is a solution.
\end{lemma}
In applying Lemma \ref{lem:zangwill} to specific algorithms, one usually needs to show the closure of the algorithm map $M$, and carefully choose the solution set $\Gamma$ and the descent function $u$. %is usually chosen as the objective function $f$.
For example, in an MM algorithm, we can choose $u$ as the objective function $f$ and the solution set
$$
    \Gamma = \{\bx\in\mathcal{X}: f(\by) \ge f(\bx), ~\forall\by\in M(\bx) \}
$$
for $M(\bx)=\argmin_{\bz\in\mathcal{X}}g(\bz|\bx)$.
Since $f(\by) \le f(\bx)$ for all $\by\in M(\bx)$ by the descent property of MM, in fact
$$
    \Gamma = \{\bx\in\mathcal{X}: f(\by) = f(\bx), ~\forall\by\in M(\bx) \} =: \mathcal{P},
$$
which we will call a set of \emph{no-progress points}.
The final requirement that $\{\bx_n\}$ is contained within a compact set is satisfied whenever $f$ is lower semicontinuous and coercive. 

We summarize the above discussion as the following proposition;
see also Proposition 8 of \cite{KeysZhouLange19ProxDist}.
\begin{proposition}[Global convergence to no-progress points]
    Suppose that the objective $f$ is lower semicontinuous and coercive, and the algorithm map $M$ defined by the MM algorithm is closed.
    Then all the limit points of the iterates $\bx_{n+1}\in M(\bx_n)$ generated by the MM algorithm are no-progress points.
\end{proposition}

This general result is slightly disappointing. Even though the objective values do not change within $\mathcal{P}$, the iterate $\{\bx_n\}$ may not even converge -- it may cycle through distinct no-progress points. % or there potentially be a point that improves $f$ outside $\mathcal{P}$.

\begin{example}[EM algorithm] As a classical example of cycling, \citet{Vaida2005} showed that in minimizing 
    $$
    f(\rho, \sigma^2) = 8\log\sigma^2 + \frac{18}{\sigma^2} + 2\log(1-\log\rho^2) + \frac{4}{\sigma^2(1-\rho^2)}
    $$
    over $\sigma^2 \le 0$ and $-1 \le \rho \le 1$ (this objective function originates from the maximum likelihood estimation of the variance and correlation coefficient of bivariate normal data with missing observations), 
    the following particular surrogate function
    $$
    g(\rho, \sigma^2\mid \rho_n, \sigma_n^2) = f(\rho, \sigma^2)
        + 2\left(\log\frac{\sigma^2(1-\rho^2)}{\sigma_n^2(1-\rho_n^2)} +  \frac{\sigma_n^2(1-\rho_n^2)}{\sigma^2(1-\rho^2)} - 1\right)
    $$
    %(recall $\log x + 1/x - 1 \ge 0$ with equality at $x=1$)
    obtained by applying the expectation-maximization (EM) algorithm, a special case of the MM algorithms,
    has two symmetric minima, $(\sigma_{n+1}^2, \rho_{n+1}) = (3, \pm\sqrt{2/3 - \sigma_n^2(1-\rho_n^2)/6})$.
    If we take $\sigma_0^2=3$ and $$\rho_{n+1} = -\sgn{(\rho_n)}\sqrt{2/3 - 3(1-\rho_n^2)/6}),$$
    then the sequence $\{(\sigma_n^2, \rho_n)\}$ oscillates between two minima $(3, \pm 1/\sqrt{3})$ of $f$ in the limit.
\end{example}
Although the above cycling can be considered desirable as it reveals multiple optima, the next example shows that this is not always the case:
\begin{example}[Generalized CCA]\label{ex:cca}
%XXX example of cycling -- Ten Berge's algorithm XXX
%As a concrete example of cycling, consider a generalization of the canonical correlation analysis (CCA).
%In its original construction, CCA	\citep{hotelling1936relations} 
%seeks directions maximizing the Pearson correlation between two sets of $n$ observations of variables of possibly different dimensions, $\bA_1\in\mathbb{R}^{n\times d_1}$ and $\bA_2\in\mathbb{R}^{n\times d_2}$, i.e.,
%\[
%	\text{maximize} ~~ 
%	(\bA_1\bt_1)^T(\bA_2\bt_2) ~~
%	\text{subject to} ~~ \|\bt_i\|=1, ~ i=1,2,
%\]
%Generalizations of CCA (i) handle more than two sets of variables $\bA_1,\dotsc,\bA_m$ (\mbox{$m\ge 2$}), and (ii) seek partial rotation matrices (as opposed to vectors) of $\bA_i$'s to achieve maximal agreement. 
The popular MAXDIFF criterion \citep{vandegeer1984linear,ten1988generalized,hanafi2006analysis} 
for generalizing the canonical correlation analysis (CCA)
into $m > 2$ sets of (partially) orthogonal matrices
solves 
\begin{equation}\label{eqn:maxdiff}
	\text{maximize} ~~ \sum_{i<j}\tr(\bO_i^T\bA_i^T\bA_j\bO_j) ~~
	\text{subject to} ~~ \bO_i^T\bO_i = \bI_r ,  ~ i=1,\dotsc,m,
\end{equation}
where $\bI_r$ is an $r\times r$ identity matrix and $\bO_i\in\mathbb{R}^{d_i\times r}$;
$\bA_i\in\mathbb{R}^{n\times d_i}$ are $n$ observations of variables of possibly different dimensions.
A standard algorithm for solving the MAXDIFF problem is Ten Berge's block relaxation algorithm \citep{ten1984orthogonal,ten1988generalized}, shown as
Algorithm \ref{alg:blockrelax}.
This is an MM algorithm (here minorization-maximization), since 
at the update of the $i$th block in the $k$th sweep, the surrogate function
\begin{align*}
	g(\bO_1,\dotsc,\bO_m |&\;\bO_1^{k+1},\dotsc,\bO_{i-1}^{k+1},\bO_i^k,\dotsc,\bO_m^k) 
	\\
	&= \frac{1}{2}\sum_{i=1}^m  \tr\left[\bO_i^T \left(\sum_{j=1}^{ i-1}\bA_{i}^T\bA_{j}\bO_j^{k+1}
	+ \sum_{j=i+1}^{ m}\bA_{i}^T\bA_{j}\bO_j^{k} \right)
	\right]
\end{align*}
minorizes the objective function of problem \eqref{eqn:maxdiff}
at 
$(\bO_1^{k+1}, \allowbreak \dotsc, \allowbreak \bO_{i-1}^{k+1}, \allowbreak \bO_i^k, \allowbreak \dotsc, \allowbreak \bO_m^k)$
and is maximized based on
the von Neuman-Fan inequality 
\[
	\tr(\bA^T\bB) \le \sum_{l} \sigma_l(\bA)\sigma_l(\bB),
\]
which holds for any two matrices $\bA$ and $\bB$ of the same dimensions with the $l$th largest singular values $\sigma_l(\bA)$ and $\sigma_l(\bB)$, respectively; equality is attained when $\bA$ and $\bB$ share a simultaneous ordered SVD \citep[see, e.g.,][]{lange2016mm}.

While each iteration monotonically improves the objective function, \citet{won2018orthogonal} show that Algorithm \ref{alg:blockrelax} may oscillate between suboptimal no-progress points. 
Set $m=3$, $d_1=d_2=d_3=d=r$ and $\bA_{1}=[\bI_d, \bI_d, \mathbf{0}]^T$, $\bA_{2}=[-\bI_d, \mathbf{0}, \bI_d]^T$,
$\bA_{3}=[\mathbf{0}, \bI_d, \bI_d]^T$
\citep[Table 1]{TenBerge1977orthogonal}.
If Algorithm \ref{alg:blockrelax} is initialized with $(\bJ,\bK,\bJ)$ where
\[
	\bJ=\begin{bmatrix} 1 & 0 \\ 0 & 1 \\ 0 & 0 \end{bmatrix}
	\quad\text{and}\quad
	\bK=\begin{bmatrix} 0 & 1 \\ 1 & 0 \\ 0 & 0 \end{bmatrix}
	,
\]
then both $\bJ-\bK$ and $\bJ+\bK$ have rank 1, and we see that $-\bK$ is one of the maximizers of $\tr\big[\bO^T(\bJ-\bK)\big]$,
and likewise
$\bJ$ maximizes $\tr\big[\bO^T(\bJ+\bK)\big]$.
Taking these values as the outputs of Line 5  of Algorithm \ref{alg:blockrelax}, we have the following cycling sequence at the end of each sweep:
\[
	(\bJ,\bK,\bJ) \to (-\bK,\bJ,-\bK) \to (-\bJ,-\bK,-\bJ) 
	\to (\bK,-\bJ,\bK) \to (\bJ,\bK,\bJ) \to \cdots.
\]
All four limit points yield the same objective values of 1. However, the global maximum of $f$ can be shown to be 3. %; see \citet{won2018orthogonal} for details.

The main reason for this oscillatory behavior is that the map $\bB=\sum_{j\neq i}\bA_{i}^T\bA_{j}\bO_j \mapsto \bP_i\bQ_i^T$ in Lines 5 and 6 is set-valued. %in general. 
If $\bB$ is rank deficient, 
any orthonormal basis of the null space of $\bB^T$ (\textit{resp.} $\bB$)  can be chosen as left (\textit{resp.} right) singular vectors corresponding to the zero singular value. Furthermore, the product $\bP_i\bQ_i^T$ may not be unique 
\citep[Proposition 7]{absil2012projection}.
\end{example}

\begin{algorithm}[h!]%[htbp]
\caption{Ten Berge's algorithm for generalized CCA}
\label{alg:blockrelax}
%\onelineup %\vspace*{-12pt}
\begin{tabbing}
{\small~~1:}\enspace Init\=ialize \= $\bO_1,\dotsc,\bO_m$\\
{\small~~2:}\enspace For $k=1,2,\dotsc$ \\
{\small~~3:}\enspace\> For $i=1,\dotsc,m$ \\
{\small~~4:}\enspace\>\> Set $\bB = \sum_{j\neq i}\bA_{i}^T\bA_j\bO_j$ \\
{\small~~5:}\enspace\>\> Compute SVD of $\bB$ as $\bP_i\bD_i\bQ_i^T$ \\
{\small~~6:}\enspace\>\> Set $\bO_i = \bP_i\bQ_i^T$ \\
{\small~~7:}\enspace\> End For \\
{\small~~8:}\enspace\> If there is no progress, then break \\
{\small~~9:}\enspace End For \\
{\small~10:}\enspace Return $(\bO_1,\dotsc,\bO_m)$
\end{tabbing}
\end{algorithm}

More satisfying ``solution sets'' are in order.
\begin{itemize}
    \item Fixed points:
    $$
        \mathcal{F} = \{\bx\in\mathcal{X}: \bx = M(\bx) \},
    $$
    if $M$ is single-valued.
    \item Stationary points:
    $$
        \mathcal{S} = \{\bx\in\mathcal{X}: d_{\bv} f(\bx) \ge 0,~\text{for all tangent vectors $\bv$ of $C$ at $\bx$}\}.
    $$
\end{itemize}
All fixed points are no-progress points, i.e., $\mathcal{F}\subset\mathcal{P}$, but not vice versa.
Note that $\by\in\mathcal{S}$ is a \emph{necessary} condition that $\by$ is a local minimizer of $f$ in $C$.
No-progress points and fixed points depend on the algorithm map $M$, whereas the stationary points depend on the problem itself.
To make $M$ single-valued, note that any convex (and weakly convex) surrogate $g(\bx|\bx_n)$ can be made strongly convex, thus attains a unique minimum, by adding the viscosity penalty $\frac{\mu}{2}\|\bx - \bx_n\|^2$ majorizing $0$; see Section \ref{sec:convergence:cycling}. If $M$ is closed and single-valued, then it is continuous. 

The classical global convergence results for MM algorithms \citep{Lange10NumAnalBook,Lange16MMBook}, which we summarize below, hinge on continuity of the map $M$:
\begin{proposition}\label{prop:fixedpoints}
    If the MM algorithm map $M$ is continuous, then $\mathcal{F}=\mathcal{P}$ and $\mathcal{F}$ is closed.
\end{proposition}
\begin{proposition}\label{prop:limitpoints}
    If (i) $f$ is continuous, (ii) $f$ is coercive, or the set $\{\bx: f(\bx) \le f(\bx_0)\}$ is compact, and (iii) the algorithm map $M$ is continuous, then every limit point of an MM sequence $\{\bx_n\}$ is a fixed point of $M$. Furthermore, $\lim_{n \to \infty} \dist(\bx_n, \mathcal{F}) = 0$.
    %and the set $\mathcal{T} \subset \mathcal{F}$ of the limit points of $\{\bx_n\}$ is closed.
\end{proposition}
\begin{proposition}\label{prop:limitpointsprops}
    Under the same assumptions of Proposition \ref{prop:limitpoints}, the MM sequence $\{\bx_n\}$ satisfies 
$$
\lim_{n \to \infty} \|\bx_{n+1} - \bx_n\| = 0.
$$
    Furthermore, the set $\mathcal{T}$ of the limit points of $\{\bx_n\}$ is compact and connected.
\end{proposition}
\noindent Note Proposition \ref{prop:limitpoints} states that $\mathcal{T}\subset\mathcal{F}$. Proposition \ref{prop:limitpointsprops} ensures there is no cycling.

Connecting the fixed points $\mathcal{F}$, which coincides with the no-progress points $\mathcal{P}$ for continuous $M$, with the stationary points $\mathcal{S}$ needs more assumptions. To equate stationary points of $f$ to those of $g(\cdot|\bx)$, we require a stronger tangency condition than the usual tangency condition \eqref{eqn:tangency}:
\begin{definition}[Strong tangency]
    An MM surrogate function $g(\cdot|\cdot)$ is said to be \emph{strongly tangent} to $f$ at $\bx\in C$ if $d_{\bv}g(\bx|\bx)=d_{\bv}f(\bx)$ for all $\bx\in C$ and tangent vector $\bv$ in $C$ at $\bx$.
\end{definition}
\begin{proposition}\label{prop:strongtangency}
    Suppose (i) the surrogate function $g(\by|\bx)$ is strongly tangent to the objective function $f$, (ii) the algorithm map $M$ is closed and single-valued, and (iii) stationary points and minimizers of $g(\by|\bx)$ are equivalent.
    Then $\mathcal{F}=\mathcal{P}=\mathcal{S}$, i.e., the sets of fixed points, no-progress points, and stationary points of $f$ coincide.
\end{proposition}
\begin{proposition}\label{prop:isolation}
    In addition to the assumptions of Proposition \ref{prop:strongtangency}, if $\mathcal{S}$, the set of all stationary points of the objective function $f$, consists of isolated points, then the set $\mathcal{T}$ of the limit points of the MM sequence $\{\bx_n\}$ is singleton, i.e., an MM sequence $\{\bx_n\}$ possesses a limit, and that limit is a stationary point of $f$ as well as a fixed point of $M$. 
\end{proposition}
\noindent Strong tangency holds when $g(\by|\bx)=f(\bx) + h(\by|\bx)$ and $h(\by|\bx)$ is differentiable with $\nabla h(\bx|\bx)=\mathbf{0}$. See \cite{Vaida2005,Yu2015} for examples of these results in action.

\iffalse
We summarize the classical global convergence results for MM algorithms \citep{Lange10NumAnalBook,Lange16MMBook}. Denote the MM algorithmic mapping $M(\bx)$.
\begin{enumerate}
\item The MM iterates $\bx_n$ monotonically decreases the objective and the values $f(\bx_n)$, if bounded below, converge. We call a point satisfying $f(M(\bx)) = f(\bx)$ a \emph{no-progress point}.

\item If (i) $f(\bx)$ is coercive, or the set $\{\bx: f(\bx) \le f(\bx_0)\}$ is compact, (ii) the algorithmic map $M(\bx)$ is continuous, and (iii) all no-progress points are fixed points of $M(\bx)$, then every cluster point of an MM sequence is a fixed point of $M(\bx)$. Furthermore, the set of fixed points $F$ is closed, and $\lim_{n \to \infty} \dist(\bx_n, F) = 0$.

\item Under same conditions (i)-(iii), the MM sequence satisfies 
$$
\lim_{n \to \infty} \|\bx_{n+1} - \bx_n\| = 0
$$
and the collection of cluster points of an MM sequence is compact and connected.

\item If all stationary points of the objective function $f(\bx)$ are isolated, then an MM sequence possesses a limit, and that limit is a stationary point of $f(\bx)$.
\end{enumerate}
\fi

We next present results that extend to non-asymptotic analysis and more general settings such as non-smooth objectives.

\subsection{Smooth objective functions}

The following proposition gives a weak form of convergence for MM
algorithms. The proposition features minimization and majorization by Lipschitz smooth functions.
\begin{proposition}
Let $f(\bx)$ be a coercive differentiable function majorized by a uniformly $L$-Lipschitz surrogate $g(\bx \mid \bx_{n})$ anchored at $\bx_{n}$. If $\by$ denotes a minimum point of $f(\bx)$, then the iterates $\bx_n$ delivered by the corresponding MM algorithm satisfy the sublinear bound
\begin{eqnarray}
\min_{0 \le k \le n}\| \nabla f(\bx_k) \|^2 & \le &
\frac{2L}{n+1} [f(\bx_0) - f(\by)] . \label{sublinear_bound}
\end{eqnarray}
When $f(\bx)$ is continuously differentiable, any limit point of the sequence $\bx_n$ is a 
stationary point of $f(\bx)$.
\end{proposition}
\begin{proof}
Given that the surrogate $g(\bx \mid \bx_n)$ satisfies the tangency condition 
$\nabla g(\bx_n \mid \bx_n) = \nabla f(\bx_n)$, the $L$-smoothness assumption entails the quadratic upper bound
\begin{eqnarray*}
f(\bx_{n+1}) - f(\bx_{n}) & \le &
g(\bx_{n+1} \mid \bx_{n}) - g(\bx_{n} \mid \bx_{n}) \\
& \le& g(\bx \mid \bx_{n}) - g(\bx_{n} \mid \bx_{n}) \\
&\le& \langle{\nabla g(\bx_n \mid \bx_{n}), \bx - \bx_{n}}\rangle
+ \frac{L}{2} \lVert{\bx - \bx_{n}}\rVert^{2} \\
&=& \langle{\nabla f(\bx_{n}), \bx - \bx_{n}}\rangle
+ \frac{L}{2} \lVert{\bx - \bx_{n}}\rVert^{2}
\end{eqnarray*}
for any $\bx$. The choice $\bx = \bx_{n} - L^{-1} \nabla f(\bx_{n})$
yields the sufficient decrease condition
\begin{eqnarray}
\label{eq:sufficient-descent}
f(\bx_{n}) - f(\bx_{n+1}) & \ge &
\frac{1}{2 L} \|\nabla f(\bx_{n})\|^{2}.
\end{eqnarray}
A simple telescoping argument now gives
\begin{eqnarray*}
\frac{n+1}{2L} \!\min_{0 \le k \le n}\| \nabla f(\bx_k) \|^2 & \!\le \!& \frac{1}{2L} \sum_{k=0}^{n} \|\nabla f(\bx_k) \|^2 \\
& \le & \amp f(\bx_{0}) - f(\bx_{n+1}) \\
& \le & \amp  f(\bx_{0}) - f(\by),
\end{eqnarray*}
which is equivalent to the bound (\ref{sublinear_bound}).
The second assertion follows directly from condition (\ref{eq:sufficient-descent}), the convergence of the sequence $f(\bx_n)$, and the continuity of $\nabla f(\bx)$.
\end{proof}

As a prelude to our next result, we state and prove a simple result of independent
interest.
\begin{proposition}
Suppose $f(\bx)$ is convex with surrogate $g(\bz \mid \bx)$ at the point $\bx$. Then $f(\bx)$ is differentiable at $\bx$ and $\nabla f(\bx)$ equals $\nabla g(\bx \mid \bx)$ wherever $\nabla g(\bx \mid \bx)$ exists.
\end{proposition}
\begin{proof} Suppose $\bx$ is such a point. Let $\bv \in \partial f(\bx)$ be a subgradient of $f(\bx)$. It suffices to show that $\bv$ 
is uniquely determined as $\bv = \nabla g(\bx \mid \bx)$. For any
direction $\bu$, consider the forward difference quotient
\begin{eqnarray*}
\frac{g(\bx+t\bu \mid \bx)-g(\bx \mid \bx)}{t}
& \ge & \frac{f(\bx+t\bu)-f(\bx)}{t} \amp \ge \amp
\langle{\bv, \bu}\rangle.
\end{eqnarray*}
Taking limits produces
$\langle{\nabla g(\bx \mid \bx), \bu}\rangle \ge 
\langle{\bv, \bu}\rangle$.
This cannot be true for all $\bu$ unless the condition
$\bv = \nabla g(\bx \mid \bx)$ holds.
\end{proof}

Imposing strong convexity on $f(\bx)$ recovers linear convergence. The ratio $\frac{\mu}{L}$ makes an appearance but, unlike convergence theorems for gradient descent, this ratio is not the condition number of either the objective or the surrogate.

\begin{proposition}
Let $f(\bx)$ be a $\mu$-strongly convex function majorized by a uniformly $L$-Lipschitz surrogate $g(\bx \mid \bx_{n})$. If the global minimum occurs at $\by$, then the MM iterates $\bx_n$ satisfy
\begin{eqnarray*}
f(\bx_{n}) - f(\by) & \le &
\left[1 - \left(\frac{\mu}{L}\right)^{2}\right]^{n}
[f(\bx_{0}) - f(\by)],
\end{eqnarray*}
thus, establishing linear convergence of $\bx_n$ to $\by$.
\end{proposition}
\begin{proof}
Existence and uniqueness of $\by$ follow from strong convexity. Because 
$\nabla g(\by \mid \by) = {\bf 0}$, the smoothness of $g(\bx \mid \by)$ gives the quadratic upper bound
\begin{eqnarray}
f(\bx) - f(\by) & \le & g(\bx \mid \by) - g(\by \mid \by) \nonumber \\
& \le & \langle{\nabla g(\by \mid \by), \bx - \by}\rangle
+ \frac{L}{2} \lVert{\bx - \by}\rVert^{2} \label{eq:prop2-pl-upper} \\
& = & \frac{L}{2} \lVert{\bx - \by}\rVert^{2}, \nonumber
\end{eqnarray}
which incidentally implies $\mu \le L$. By the previous proposition, $f(\bx)$ is everywhere differentiable with $\nabla f(\bx)=
\nabla g(\bx \mid \bx)$. In view of the strong convexity assumption, we have the lower bound
\begin{eqnarray}
\lVert{\nabla f(\bx)}\rVert \cdot \lVert{\by - \bx}\rVert
&\ge & - \langle{\nabla f(\bx), \by - \bx}\rangle \nonumber \\
&\ge &  f(\by) - f(\bx) - \langle{\nabla f(\bx), \by - \bx}\rangle
\label{eq:prop2-pl-lower} \\
&\ge & \frac{\mu}{2} \lVert{\by - \bx}\rVert^{2}. \nonumber
\end{eqnarray}
It follows that $\lVert{\nabla f(\bx)}\rVert \ge \frac{\mu}{2}\lVert{\by - \bx}\rVert$. Combining inequalities (\ref{eq:prop2-pl-upper}) and (\ref{eq:prop2-pl-lower}) furnishes the Polyak-\L{}ojasiewicz (PL) bound
\begin{eqnarray*}
\lVert{\nabla f(\bx)}\rVert^{2} & \ge &
\frac{\mu^{2}}{2 L} [f(\bx) - f(\by)].
\end{eqnarray*}
We now turn to the MM iterates and take $\bx = \bx_{n} - \frac{1}{L} \nabla f(\bx_{n})$. The PL inequality implies
\begin{eqnarray*}
f(\bx_{n+1}) - f(\bx_{n})
& \le & g(\bx_{n+1} \mid \bx_{n}) - g(\bx_{n} \mid \bx_{n}) \\
& \le & g(\bx \mid \bx_{n}) - g(\bx_{n} \mid \bx_{n}) \\
& \le & \Big\langle{\nabla g(\bx_{n} \mid \bx_{n}), -\frac{1}{L} \nabla g(\bx_{n} \mid \bx_{n})}\Big\rangle +
\frac{L}{2}\Big\lVert{\frac{1}{L} \nabla g(\bx_{n} \mid \bx_{n})}\Big\rVert^{2} \\
& = & -\frac{1}{2L}\lVert{\nabla f(\bx_{n})}\rVert^{2} \\
&\le & - \frac{2\mu^2}{L^2} \left[ f(\bx_{n}) - f(\by)\right].
\end{eqnarray*}
Subtracting $f(\by)$ from both sides of the previous inequality 
and rearranging gives
\begin{eqnarray*}
f(\bx_{n+1}) - f(\by) & \le & \left[1 - \frac{\mu^{2}}{2 L^{2}}
    \right]
    [f(\bx_{n}) - f(\by)].
\end{eqnarray*}
Iteration of this inequality yields the claimed linear convergence.
\end{proof}

%$K_1^o$ = metric inequality constraints, $K_2^o = \mathbb{R}^p_+$.
%Since $(K_1+K_2)^o = K_1^0 \cap K_2^0$ it suffices to project onto
%$(K_1+K_2)^o$. But then it suffices to project onto $K_1+K_2$. This 
%can be done by the Minkowski alternating projections.

\subsection{Non-smooth objective functions}

Consider an MM minimization algorithm with objective $f(\bx)$ and surrogates $g(\bx \mid \bx_n)$. If $f(\bx)$ is coercive and continuous and the $g(\bx \mid \bx_n)$ are $\mu$-strongly convex, then we know that the MM iterates $\bx_{n+1} = \argmin_{\bx} g(\bx \mid \bx_n)$ remain within the compact sublevel set $S = \{\bx: f(\bx) \le f(\bx_0)\}$ \citep{Lange16MMBook}. Furthermore, the strong convexity inequality
\begin{eqnarray}
f(\bx_n)-f(\bx_{n+1}) & \ge & g(\bx_n \mid \bx_n) - g(\bx_{n+1}\mid \bx_n)
\amp \ge \amp \frac{\mu}{2}\|\bx_n-\bx_{n+1}\|^2
\label{strong_convexity_ineq}
\end{eqnarray}
implies that 
\begin{eqnarray*}
\sum_{n=0}^\infty \|\bx_n-\bx_{n+1}\|^2
& \le & \frac{1}{\mu} \Big[f(\bx_0)-\bar{f}\Big],
\end{eqnarray*}
where $\bar{f} = \lim_{n \to \infty} f(\bx_n)$. It follows from a well-known theorem of Ostrowski \citep{Lange16MMBook} that the set $W$ of limit points of the $\bx_n$ is compact and connected. It is also easy to show that $f(\bx)$ takes the constant value $\bar{f}$ on $W$ and that 
$\lim_{n \to \infty} \dist(\bx_n, W) = 0$.

We will need the concept of a Fr\'echet subdifferential. If $f(\bx)$
is a function mapping $\mathbb{R}^p$ into $\mathbb{R} \cup \{+\infty\}$,
then its Fr\'echet subdifferential at $\bx \in \dom f$ is the set
\begin{eqnarray*}
\partial^F f(\bx) & = & \left\{\bv: \liminf_{\by \to \bx}
\frac{f(\by)-f(\bx) - \bv^t(\by-\bx)}{\|\by-\bx\|} \ge 0 \right\}.
\end{eqnarray*}
The set $\partial^F f(\bx)$ is closed, convex, and possibly empty. If $f(\bx)$ is convex, then $\partial^F f(\bx)$ reduces to its convex subdifferential. If $f(\bx)$ is differentiable, then $\partial^F f(\bx)$ reduces to its ordinary differential. At a local minimum $\bx$, Fermat's rule ${\bf 0} \in \partial^F f(\bx)$ holds. 

\begin{proposition}
In an MM algorithm, suppose $f(\bx)$ is coercive, $g(\bx \mid \bx_n)$ is differentiable, and the algorithm map $M(\bx)$ is closed. Then all points 
$\bz$ of the convergence set $W$ are critical in the sense that 
${\bf 0} \in \partial^F (-f)(\bz)$. 
\end{proposition}
\begin{proof} 
Let the subsequence $x_{n_m}$ of the MM sequence 
$\bx_{n+1} \in M(\bx_n)$ converge to $\bz \in W$. By passing a subsubsequence if necessary, we may suppose that $\bx_{n_m+1}$ converges to $\by$. Owing to our closedness assumption, $\by \in M(\bz)$. Given that $f(\by)=f(\bz)$, it is obvious that $\bz$ also minimizes $g(\bx \mid \bz)$ and that 
${\bf 0} = \nabla g(\bz \mid \bz)$.  Since the difference $h(\bx \mid \bz)= g(\bx \mid \bz)-f(\bx)$ achieves it minimum at $\bx = \bz$, the Fr\'echet subdifferential $\partial^F h(\bx \mid \bz)$ satisfies 
\begin{eqnarray*}
{\bf 0} & \in & \partial^F h(\bz \mid \bz) \amp = \amp
\nabla g(\bz \mid \bz) + \partial^F (-f) (\bz).
\end{eqnarray*}
It follows that ${\bf 0} \in \partial^F (-f) (\bz)$.
\end{proof}

We will also need to invoke \L{}ojasiewicz's inequality. This deep result depends on some rather arcane algebraic geometry \citep{BierstoneMilman88Semianalytic,BochnakCosteRoy98RealAlgebraicGeometry}. It applies to semialgebraic functions and their more inclusive cousins semianalytic functions and subanalytic functions. For simplicity we focus on semialgebraic functions. The class of semialgebraic subsets of $\mathbb{R}^p$ is the smallest class such that:
\begin{description}
\item[a)] It contains all sets of the form
$\{\bx: q(\bx)>0 \}$ for a polynomial $q(\bx)$ in $p$ variables.
\item[b)] It is closed under the formation of finite unions, finite intersections, and set complementation. 
\end{description}
A function $a:\mathbb{R}^p \mapsto \mathbb{R}^r$ is said to be semialgebraic if its graph is a semialgebraic set of $\mathbb{R}^{p+r}$.
The class of real-valued semialgebraic contains all polynomials $p(\bx)$.
It is closed under the formation of sums, products, absolute values,
reciprocals when $a(\bx) \ne 0$, roots when $a(\bx) \ge 0$, and maxima 
$\max\{a(\bx),b(\bx)\}$ and minima $\min\{a(\bx),b(\bx)\}$. For our purposes, it is important to note that $\dist(\bx, S)$ is a semialgebraic function whenever $S$ is a semialgebraic set.

\L{}ojasiewicz's inequality in its modern form \citep{AttouchBolte09KLIneq} requires that $f(\bx)$ be continuous and subanalytic with a closed domain. If $\bz$ is a critical point of $f(\bx)$, then 
\begin{eqnarray*}
|f(\bx) - f(\bz)|^{\theta(\bz)} & \le & c(\bz) \|\bv \|
\end{eqnarray*}
for some constant $c(\bz)$, all $\bx$ in some open ball $B_{r(\bz)}(\bz)$ around $\bz$ of radius $r(\bz)$, and all $\bv$ in $\partial^F f(\bx)$. This inequality applies to semialgebraic functions since they are automatically subanalytic. We will apply \L{}ojasiewicz's inequality to the points in the limit set $W$.

\subsubsection{MM convergence for semialgebraic functions}

\begin{proposition} Suppose that $f(\bx)$ is coercive, continuous, and subanalytic and all $g(\bx \mid \bx_n)$ are continuous, $\mu$-strongly convex, and satisfy the
Lipschitz condition
\begin{eqnarray*}
\|\nabla g(\bu \mid \bx_n)-\nabla g(\bv \mid \bx_n) \|
& \le & L \|\bu-\bv\| 
\end{eqnarray*}
on the compact sublevel set $\{\bx: f(\bx) \le f(\bx_0)\}$. Then the MM iterates 
$\bx_{n+1} = \argmin_{\bx} g(\bx \mid \bx_n)$ converge to a critical point in $W$.  
\end{proposition} 
\begin{proof}
Because $h(\bx \mid \by)= g(\bx \mid \by)-f(\bx)$ achieves it minimum at $\bx = \by$, the Fr\'echet subdifferential $\partial^F h(\bx \mid \by)$ satisfies 
\begin{eqnarray*}
{\bf 0} & \in & \partial^F h(\by \mid \by) \amp = \amp
\nabla g(\by \mid \by) + \partial^F (-f) (\by).
\end{eqnarray*}
It follows that $-\nabla g(\by \mid \by) \in \partial^F (-f) (\by)$.
By assumption
\begin{eqnarray*}
\|\nabla g(\bu \mid \bx_n)-\nabla g(\bv \mid \bx_n) \|
& \le & L\|\bu-\bv\| 
\end{eqnarray*}
for all $\bu$ and $\bv$ and $\bx_n$. In particular, because $\nabla g(\bx_{n+1} \mid \bx_n) = {\bf 0}$, we have
\begin{eqnarray}
\|\nabla g(\bx_n \mid \bx_n) \|
& \le & L \|\bx_{n+1}-\bx_n\|. \label{lipschitz_MM}
\end{eqnarray}

According to the \L{}ojasiewicz inequality applied for the subanalytic function $\bar{f}-f(\bx)$, for each $\bz \in W$ there exists a radius $r(\bz)$ and an exponent $\theta(\bz) \in [0,1)$ with 
\begin{eqnarray*}
|f(\bu) - f(\bz)|^{\theta(\bz)} & = & 
|\bar{f} - f(\bu) - \bar{f}+\bar{f}|^{\theta(\bz)}
\amp \le \amp c(\bz) \|\bv \|
\end{eqnarray*}
for all $\bu$ in the open ball $B_{r(\bz)}(\bz)$ around $\bz$ of radius $r(\bz)$ and all 
$\bv \in \partial^F (\bar{f}-f) (\bu)= \partial^F (-f) (\bu)$. We will apply this inequality to $\bu=\bx_n$ and 
$\bv = -\nabla g(\bx_n \mid \bx_n)$.
In so doing, we would like to assume that the exponent $\theta(\bz)$ and constant $c(\bz)$ do not depend on $\bz$. With this end in mind, cover 
$W$ by a finite number of balls balls $B_{r(\bz_i)}(\bz_i)$ and take 
$\theta=\max_i \theta(\bz_i)<1$ and $c = \max_i c(\bz_i)$. For a sufficiently large $N$, every $\bx_n$ with $ n \ge N$ falls within one
of these balls and satisfies $|\bar{f}-f(\bx_n)| < 1$. Without loss
of generality assume $N=0$. The \L{}ojasiewicz inequality now
entails
\begin{eqnarray}
|\bar{f}-f(\bx_n)|^{\theta} & \le & c \|\nabla g(\bx_n \mid \bx_n) \|. 
\label{Lojasiewicz_MM}
\end{eqnarray}
In combination with the concavity of the function $t^{1-\theta}$ on $[0,\infty)$, inequalities (\ref{strong_convexity_ineq}),
(\ref{lipschitz_MM}), and (\ref{Lojasiewicz_MM}) imply
\begin{eqnarray*}
[f(\bx_n) - \bar{f}]^{1-\theta}-[f(\bx_{n+1})-\bar{f}]^{1-\theta} & \ge & 
\frac{1-\theta}{[f(\bx_n)-\bar{f}]^{\theta}}[f(\bx_n)-f(\bx_{n+1})] \\
& \ge & \frac{1-\theta}{c\|\nabla g(\bx_n \mid \bx_n) \|}
\frac{\mu}{2}\|\bx_{n+1}-\bx_n\|^2 \\
& \ge & \frac{(1-\theta)\mu}{2cL}\|\bx_{n+1}-\bx_n\|.
\end{eqnarray*}
Rearranging this inequality and summing over $n$ yield
\begin{eqnarray*}
\sum_{n=0}^\infty \|\bx_{n+1}-\bx_n\| & \le &
\frac{2cL}{(1-\theta)\mu}[f(\bx_0) - \bar{f}]^{1-\theta}
\end{eqnarray*}
Thus, the sequence $\bx_n$ is a fast Cauchy sequence and converges to a unique limit in $W$.
\end{proof}

\subsection{A proximal trick to prevent cycling}\label{sec:convergence:cycling}

Consider minimizing a function $f(\bx)$ bounded below and possibly subject to constraints. The MM principle involves constructing a surrogate function $g(\bx \mid \bx_n)$ that majorizes $f(\bx)$ around $\bx_n$. For any $\rho>0$, adding the penalty $(\rho/2) \|\bx - \bx_n\|^2$ to the surrogate produces a new surrogate
\begin{eqnarray*}
g(\bx \mid \bx_n) + \frac{\rho}{2} \|\bx - \bx_n\|^2.
\end{eqnarray*}
Rearranging the inequality
\begin{eqnarray*}
	g(\bx_{n+1} \mid \bx_n) + \frac{\rho}{2} \|\bx_{n+1} - \bx_n\|^2 \le g(\bx_n \mid \bx_n)
\end{eqnarray*}
yields
\begin{eqnarray*}
	\frac{\rho}{2}  \|\bx_{n+1} - \bx_n\|^2 \le g(\bx_n \mid \bx_n) - g(\bx_{n+1} \mid \bx_n) \le f(\bx_n) - f(\bx_{n+1}).
\end{eqnarray*}
Thus the MM iterates induced by the new surrogate satisfy
\begin{eqnarray*}
\lim_{n \to \infty} \|\bx_{n+1} - \bx_n\| = 0.
\end{eqnarray*}
This property is inconsistent with algorithm cycling between distant limit points.

\section{Paracontraction}

Another useful tool for proving iterate convergence of MM algorithms is paracontraction. Recall that a map $T : \mathcal{X} \to \mathbb{R}^d$ is contractive with respect to a norm $\|\bx\|$ if $\|T(\by) - T(\bz)\| < \|\by - \bz\|$ for all $\by \neq \bz$ in $\mathcal{X}$. It is strictly contractive if there exists a constant $c \in [0,1)$ with $\|T(\by) - T (\bz)\| \le c\|\by - \bz\|$ for all such pairs. If $c=1$, then the map is nonexpansive.
\begin{definition}[Paracontractive map]
    A map $T : \mathcal{X} \to \mathbb{R}^d$ is said to be \emph{paracontractive} if for every fixed point $\by$ of $T$ (i.e., $\by=T(\by)$), 
    the inequality $\|T(\bx) - \by\| < \|\bx - \by\|$ holds unless $\bx$ is itself a fixed point. 
\end{definition}
\noindent A strictly contractive map is contractive, and a contractive map is paracontractive.

An important result regarding paracontractive maps is the theorem of Elsner, Koltract, and Neumann \citep{elsner92}, which states that whenever a continuous paracontractive map $T$ possesses one or more fixed points, then the sequence of iterates $\bx_{n+1} = T(\bx_n)$ converges to a fixed point regardless of the initial point $\bx_0$. More formal statement is as follows:
\begin{proposition}\label{prop:elsner}
Suppose the continuous maps $T_0, \cdots, T_{r-1}$ of a set into itself are paracontractive under the norm $\|x\|$. 
Let $F_i$ denote the set of fixed points of $T_i$. 
If the intersection $F = \cap_{i=0}^{r-1}F_i$ is nonempty, then the sequence
$$
    \bx_{n+1} = T_{n\mod r}(\bx_n)
$$
converges to a limit in $F$. 
In particular, if $r = 1$ and $T = T_0$ has a nonempty set of fixed points $F$, then 
$\bx_{n+1} = T(\bx_n)$ converges to a point in $F$.
\end{proposition}
\noindent A simple proof is given in \citet{lange13}.

Proposition \ref{prop:elsner} converts the task of proving convergence of MM iterates to that of showing 1) continuity, 2) paracontractivity, and 3) existence of a fixed point, of the MM algorithm map,
and that 4) any fixed point is a stationary point of the objective.
A nice example is the recent work by \citet{won2019projection} on Euclidean projection onto the Minkowski sum of sets. The Minkowski sum of two sets $A$ and $B$ in $\mathbb{R}^d$ is
$$
A+B = \{\ba+\bb: \ba \in A, ~ \bb \in B\}.
$$
It is easy to show that $A+B$ is convex whenever $A$ and $B$ are both convex and is closed if at least one of the two sets is compact and the other is closed. 
When $A+B$ is closed with $A$ and $B$  convex, we may employ a block descent algorithm, an instance of MM algorithms, for finding the closest point to $\bx \notin A+B$, which consists of alternating 
\begin{equation}\label{eqn:blockdescent}
\begin{split}
	\bb_{n+1} & =  P_B(\bx - \ba_n) \\
	\ba_{n+1} & =  P_A(\bx - \bb_{n+1}),
\end{split}
\end{equation}
assuming that the projection operators $P_A$ and $P_B$ onto $A$ and $B$ are both known or easy to compute.

In order to show that the sequence
$\{\ba_n+\bb_n\}$ converges to the closest point using Proposition \ref{prop:elsner}, we first need to show the continuity of the map
$$
T(\ba) =  P_A[\bx - P_B(\bx-\ba)].
$$
The obtuse angle property of Euclidean projection \citep[Example 6.5.3]{lange13} yields
\begin{align*}
    \langle \ba - P_A(\ba), P_A(\tilde{\ba}) - P_A(\ba) \rangle &\le 0 \\
    \langle \tilde{\ba} - P_A(\tilde{\ba}), P_A(\ba) - P_A(\tilde{\ba}) \rangle &\le 0 
\end{align*}
for any $\ba, \tilde{\ba} \in \mathbb{R}^d$.
Adding these inequalities, rearranging, and applying the Cauchy-Schwarz inequality give
\begin{equation}\label{eqn:nonexpansive}
\begin{split}
    \|P_A(\ba) - P_A(\tilde{\ba})\|^2
    &\le \langle \ba - \tilde{\ba}, P_A(\ba) - P_A(\tilde{\ba}) \rangle
    \\
    &\le \|\ba - \tilde{\ba}\|\|P_A(\ba) - P_A(\tilde{\ba})\|
    .
\end{split}
\end{equation}
Thus $\|P_A(\ba) - P_A(\tilde{\ba})\| \le \|\ba - \tilde{\ba}\|$. That is, $P_A$ is nonexpansive, and the inequalty holds if and only if
\begin{equation}\label{eqn:cauchyswartz}
    P_A(\ba) - P_A(\tilde{\ba}) = c(\ba - \tilde{\ba})
\end{equation}
for some constant $c$. Likewise, $P_B$ is nonexpansive.
%
%Because both $P_A$ and $P_B$ are nonexpansive with respect to the Euclidean norm,
Therefore,
\begin{equation}\label{eqn:paracontraction}
\begin{split}
& \|P_A[\bx - P_B(\bx-\ba)] - P_A[\bx - P_B(\bx-\tilde{\ba})]\|
\\
& \le 
\|P_B(\bx-\ba)- P_B(\bx-\tilde{\ba})\|
\le 
\|\ba-\tilde{\ba}\|.
\end{split}
\end{equation} 
This proves that $T$ is nonexpansive, hence continuous.

Next, we show that $T$ is paracontractive. Suppose $\tilde{\ba}$
is a fixed point, $\ba \neq \tilde{\ba}$, and equality holds throughout inequalities \eqref{eqn:paracontraction}.
%The standard proof that a convex set projection is paracontractive \citep[pp. 389--390]{lange13}
Inequalities \eqref{eqn:nonexpansive} and equation \eqref{eqn:cauchyswartz}
indicate that equality is achieved in the previous two inequalities only if
\begin{align*}
&P_A[\bx - P_B(\bx-\ba)]-[\bx - P_B(\bx-\ba)]
\\
& ~= 
P_A[\bx - P_B(\bx-\tilde{\ba})]-[\bx - P_B(\bx-\tilde{\ba})], 
\end{align*}
and
\begin{align*}
P_B(\bx-\ba)-(\bx-\ba) 
 = P_B(\bx-\tilde{\ba})-(\bx-\tilde{\ba}).
\end{align*} 
Subtracting the second of these equalities from the first
gives
\[
P_A[\bx - P_B(\bx-\ba)]-\ba 
= 
P_A[\bx - P_B(\bx-\tilde{\ba})]-\tilde{\ba} 
= {\bf 0}.
\]
It follows that equality in inequalities \eqref{eqn:paracontraction} is achieved only if $\ba$ is also
a fixed point. 

To show that $T$ possesses a fixed point,
note that given the closedness of $A+B$, there exists a closest point
$\tilde{\ba}+\tilde{\bb}$ to $\bx$, where $\tilde{\ba}\in A$ and $\tilde{\bb}\in B$.
Since block descent cannot improve
the objective $f(\ba,\bb) = \frac{1}{2}\|\bx-\ba-\bb\|^2$ on the set $A\times B$ starting from $(\tilde{\ba},\tilde{\bb})$, it is clear that $\tilde{\ba} = T(\tilde{\ba})$.

Finally, suppose $\tilde{\ba}$ is any fixed point, and define $\tilde{\bb} = P_B(\bx-\tilde{\ba})$.
To prove that $\tilde{\ba}+\tilde{\bb}$ minimizes the distance to $\bx$, it suffices
to show that for every tangent vector $\bv =\ba+\bb-\tilde{\ba}-\tilde{\bb}$ at
$\tilde{\ba}+\tilde{\bb}$, the directional derivative
\begin{align*}
d_{\bv}\frac{1}{2} \|\bx-\tilde{\ba}-\tilde{\bb}\|^2
& =  -\langle \bx-\tilde{\ba}-\tilde{\bb}, \bv \rangle \\
& =  -\langle \bx-\tilde{\ba}-\tilde{\bb}, \ba-\tilde{\ba} \rangle - \langle \bx-\tilde{\ba}-\tilde{\bb}, \bb-\tilde{\bb} \rangle
\end{align*} 
is nonnegative. However, the inequalities 
$-\langle \bx-\tilde{\ba}-\tilde{\bb}, \ba-\tilde{\ba}\rangle \ge 0$
and $-\langle\bx-\tilde{\ba}-\tilde{\bb}, \bb-\tilde{\bb}\rangle \ge 0$ hold because $\tilde{\ba}$ minimizes
$\ba \mapsto \frac{1}{2}\|\bx-\ba-\tilde{\bb}\|^2$ and $\tilde{\bb}$ minimizes
$\bb \mapsto \frac{1}{2}\|\bx-\tilde{\ba}-\bb\|^2$.
Thus, any fixed point of $T$ furnishes a minimum of the convex function
$f(\ba,\bb)$ on the set $A \times B$.

\section{Bregman Majorization}

Bregman majorization is a technique for constructing a sequence of surrogate functions pertinent to an MM algorithm. Let us first define the notion of Bregman divergence. 

\begin{definition}[Bregman divergence]
	For a proper convex function $\phi(\bx)$ that is continuously differentiable on $\mathbf{int}\dom{\phi}$, the Bregman divergence $B_{\phi}:\mathcal{X}\times\mathcal{X}\to\mathbb{R}$ is defined as
	$$
	B_{\phi}(\bx\Vert\by) = \phi(\bx) - \phi(\by) - \langle \nabla\phi(\by), \bx - \by \rangle, 
	\quad \bx, \by \in \mathbf{int}\dom{\phi}
	.
	$$
\end{definition}

We are concerned with the following optimization problem:
\begin{equation}\label{eqn:main}
	\min_{\bx\in C} f(\bx), \quad C \subset \mathcal{X} \text{ is closed and convex,}
\end{equation}
where $f(\bx)$ is a convex, proper, and lower semicontinuous.
In order to solve this problem, the Bregman majorization method constructs the sequence of surrogate functions 
$$
	g(\bx\mid\bx_n) = f(\bx) + B_{\phi}(\bx\Vert\bx_n)
	.
$$
and successively minimizes these.
This is a valid MM algorithm since
the following properties of the Bregman divergence are immediate from definition:
\begin{enumerate}
	\item $B_{\phi}(\bx\Vert\by) \ge 0$;
	\item $B_{\phi}(\bx\Vert\bx) = 0$;
	\item If $\phi$ is strictly convex, then $B_{\phi}(\bx\Vert\by)=0$ if and only if $\bx=\by$.
\end{enumerate}
Thus, $g(\bx\mid\bx_n) \ge f(\bx)$ for all $\bx$ and $g(\bx_n\mid\bx_n) = f(\bx_n)$.
We can choose $\phi(\bx)$ so that for $\mathbf{cl}\dom{\phi} = C$.

The subsequent section studies the convergence property of the Bregman majorization. 
\subsection{Convergence analysis via SUMMA}

The sequential unconstrained minimization method algorithm \citep[SUMMA;][]{Byrne08Sequential} is a class of algorithms for solving  optimization problems of the form
\begin{equation}\label{eqn:constrained}
	\min_{\bx\in C} f(\bx), \quad C \subset \mathcal{X} \text{ is closed,}
\end{equation}
by minimizing a sequence of auxiliary functions 
$$
	G_n(\bx) = f(\bx) + g_{n+1}(\bx), \quad n=1, 2, \dotsc,
$$
over $\mathcal{X}$. 
The minimizer of $G_n(\bx)$ is denoted by $\bx_n$.
The conditions imposed on the sequence of functions $g_n(\bx)$ are: 
\begin{enumerate}
	\item $g_n(\bx) \ge 0$ for all $\bx \in \mathcal{X}$;
	\item $g_n(\bx_{n-1}) = 0$; 
	\item $G_n(\bx) - G_n(\bx_n) \ge g_{n+1}(\bx)$ for all $\bx \in C$.
\end{enumerate}

If $g_n(\bx)$ depends on $n$ only through the iterate $x_n$, then this method coincides with the MM algorithm by identifying $G_n(\bx) = g(\bx \mid \bx_n)$ and $g_n(\bx) = g(\bx \mid \bx_{n-1}) - f(\bx)$, with the additional requirement
\begin{equation}\label{eqn:summa}\tag{SUMMA}
	g(\bx \mid \bx_n) - g(\bx_{n+1} \mid \bx_n) \ge g(\bx \mid \bx_{n+1}) - f(\bx)
\end{equation}
for all $\bx \in C$.

Let us show that condition \eqref{eqn:summa} is satisfied by the Bregman majorization
%. To see this, note that $\bx_{n+1}$ minimizes 
$g(\bx \mid \bx_n) = 
%f(\bx) + B_{\phi}(\bx \mid %\bx_n) = 
f(\bx) + \phi(\bx) - \phi(\bx_n) - \langle \nabla\phi(\bx_n), \bx - \bx_n \rangle$. The optimality condition for minimizing $g(\bx \mid \bx_n)$ is
$$
	{\bf 0} \in \partial f(\bx_{n+1}) + \nabla\phi(\bx_{n+1}) - \nabla\phi(\bx_n).
$$
For the appropriate choice of $\bs_{n+1} \in \partial f(\bx_{n+1})$, it follows that
\begin{align*}
	g(\bx \mid \bx_n) - g(\bx_{n+1} \mid \bx_n)
   %	&= f(\bx) + B_{\phi}(\bx\Vert\bx_n) - f(\bx_{n+1}) - B_{\phi}(\bx_{n+1}\Vert\bx_n) \\
%	&= f(\bx) - f(\bx_{n+1})  \\
%	& \quad - \phi(\bx_{n+1}) + \phi(\bx_n) + \langle \nabla\phi(\bx_n), \bx_{n+1} - %\bx_n \rangle  \\
%	& \quad + \phi(\bx) - \phi(\bx_n) - \langle \nabla\phi(\bx_n), \bx - \bx_n %\rangle \\
	&= f(\bx) - f(\bx_{n+1})  + \phi(\bx) - \phi(\bx_{n+1}) \\
	& - \langle  \nabla\phi(\bx_{n}), \bx - \bx_{n+1} \rangle \\
	&= f(\bx) - f(\bx_{n+1})  - \langle \bs_{n+1}, \bx - \bx_{n+1} \rangle \\
	& \quad + \phi(\bx) - \phi(\bx_{n+1}) - \langle \nabla\phi(\bx_{n+1}), \bx - \bx_{n+1} \rangle \\
	&\ge B_{\phi}(\bx\Vert\bx_{n+1}) = g(\bx\mid\bx_{n+1}) - f(\bx),
\end{align*}
where the last inequality is a consequence of the convexity of $f(\bx)$.

The following propositions concern convergence of MM algorithms satisfying condition \eqref{eqn:summa}.
\begin{proposition}\label{prop:objectiveconvergence}
	Assume (a) $p^\star = \inf_{x\in C}f(\bx) > -\infty$ and 
    (b) $\bx_n \in C$ for each $n$.
	If condition \eqref{eqn:summa} holds, then any MM sequence generated by the map $\bx_{n+1} \in \argmin_{\bx\in\mathcal{X}} g(\bx \mid \bx_n)$ satisfies $\lim_{n \to \infty}f(\bx_n) = p^\star$.
\end{proposition}
\begin{proof}
	By the descent property of MM and the bound $f(\bx_n) \ge p^\star > -\infty$ given $\bx_n \in C$, 
	the sequence $f(\bx_n)$ converges to a limit $d \ge p^\star$.
	Suppose for some $\bx\in C$ that $f(\bx) < d$. Then, by condition \eqref{eqn:summa},
	\begin{align*}
		[g(\bx \mid \bx_n) - f(\bx)] - [g(\bx \mid \bx_{n+1}) - f(\bx)] &\ge g(\bx_{n+1} \mid \bx_n) - f(\bx) \\
			&\ge f(\bx_{n+1}) - f(\bx) \\
			&\ge d - f(\bx) \\
			& > 0.
	\end{align*}
	Thus. the sequence $g(\bx \mid \bx_n) - f(\bx)$ decreases and its successive differences are bounded away from zero. The latter property contradicts the requirement for the surrogate function that $g(\bx \mid \bx_n) \ge f(\bx)$, and therefore $d=p^\star$.
	%Therefore $f(\bx) \ge d$ for all $\bx \in C$ or $p^\star = \inf_{\bx\in C}f(\bx) \ge d \ge p^\star$. That is, 
\end{proof}
\begin{proposition}\label{prop:iterateconvergence}
	In addition to the assumptions of Proposition \ref{prop:objectiveconvergence},
	further assume that (c) the minimum $p^\star$ is attained and the set $F$ of the minimizers of $f(\bx)$ in $C$ is nonempty,	(d) $f(\bx)$ is continuous on $D\subset\mathcal{X}$ such that $\mathbf{cl}D = C$, (e) for each $n$ $g(\bx \mid \bx_n)$ is $\mu$-strongly convex with respect to the norm $\|\cdot\|$ and $\dom g(\cdot \mid \bx_n)=D$, and 
	(f) $g(\bx \mid \bx_n) - f(\bx) \le \frac{L}{2}\|\bx - \bx_n\|^2$ for all $\bx \in D$ and each $n$.
	If condition \eqref{eqn:summa} holds, then the MM sequence $\bx_{n+1} = \argmin_{\bx\in\mathcal{X}} g(\bx \mid \bx_n)$ converges to a point in $F$.
\end{proposition}
\begin{proof}
	Because of strong convexity, the minimum of $g(\bx \mid \bx_n)$ is uniquely attained for each $n$. Furthermore, for any $\bx \in D$,
	\begin{equation}\label{eqn:strongconvexity}
	g(\bx \mid \bx_n) - g(\bx_{n+1} \mid \bx_n) \ge \frac{\mu}{2}\|\bx - \bx_{n+1}\|^2
	\end{equation}
	Let $\by \in F$ be a minimizer of $f(\bx)$ in $C$. 
	Since $f(\bx_{n+1}) \le g(\bx_{n+1} \mid \bx_n)$,
	\begin{equation}\label{eqn:lowerbound}
		g(\by \mid \bx_n) - f(\bx_{n+1}) \ge g(\by \mid \bx_n) - g(\bx_{n+1} \mid \bx_n) 
			\ge \frac{\mu}{2}\|\by - \bx_{n+1} \|^2,
	\end{equation}
	where the last inequality follows from the strong convexity  of $g(\bx \mid \bx_n)$.
	Condition \eqref{eqn:summa} also implies
	\begin{align*}
		[g(\by \mid \bx_n) - f(\by)] - [g(\by \mid \bx_{n+1}) - f(\by)] &\ge g(\bx_{n+1} \mid \bx_n) - f(\by) \\
			&\ge f(\bx_{n+1}) - p^\star 
			\ge 0
	\end{align*}
	Hence the decreasing nonnegative sequence $g(\by \mid \bx_n) - f(\by)$ has a limit. 	In addition, $f(\by) - f(\bx_{n+1})$ tends to zero by Proposition \ref{prop:objectiveconvergence}.
	It follows that the leftmost side of inequality \eqref{eqn:lowerbound} tends to a limit, and the sequence $\bx_n$ is bounded.

	Suppose the convergent subsequence $\bx_{n_m}$ of $\bx_n$ has a limit $\bz$. By continuity, $f(\bz) = \lim_{m\to\infty} f(\bx_{n_m}) = p^\star$, so $\bz$ is also optimal.
	Now,
	\begin{align*}
		0 &\le g(\bz \mid \bx_n) - g(\bx_{n+1} \mid \bx_n) \\
			&= [g(\bz \mid \bx_n) - f(\bz)] + f(\bz) - f(\bx_{n+1}) - [g(\bx_{n+1} \mid \bx_n) - f(\bx_{n+1})] \\
			&\le g(\bz \mid \bx_n) - f(\bz) \\
			&\le \frac{L}{2}\|\bx_n - z\|^2
	\end{align*}
	due to $f(\bz) \le f(\bx_{n+1})$, $g(\bx_{n+1} \mid \bx_n) - f(\bx_{n+1}) \ge 0$, and assumption (f).
	Again by Condition \eqref{eqn:summa}, we further have
	\begin{equation}\label{eqn:sandwich}
		0 \le g(\bz \mid \bx_n) - g(\bx_{n+1} \mid \bx_n) \le g(\bz \mid \bx_n) - f(\bz) 
		\le g(\bz \mid \bx_{n-1}) - g(\bx_{n} \mid \bx_{n-1}). 
	\end{equation}
	Thus, the nonnegative sequence $g(\bz \mid \bx_n) - g(\bx_{n+1} \mid \bx_n)$ is monotonically decreasing and convergent. Its subsequence $g(\bz \mid \bx_{n_m}) - g(\bx_{n_m+1} \mid \bx_{n_m})$ is also bounded by $\frac{L}{2}\|\bx_{n_m} - \bz\|^2$, which converges to zero. Thus the whole sequence tends to zero. 
	By inequality \eqref{eqn:sandwich}, it follows that the sequence $g(\bz \mid \bx_n) - f(\bz)$ converges to zero.

	The final inequality
	\begin{align*}
		g(\bz \mid \bx_n) - f(\bz) &= g(\bz \mid \bx_n) - g(\bx_{n+1} \mid \bx_n) + g(\bx_{n+1} \mid \bx_{n}) - f(\bz) \\
		&\ge \frac{\mu}{2}\|\bz - \bx_{n+1}\|^2 + f(\bx_{n+1}) - f(\bz)
	\end{align*}
	now proves that the entire sequence $\bx_n$ converges to $\bz \in F$.
\end{proof}

\begin{remark}
	Assumption (e) (uniform strong convexity of the surrogate functions) is much less restrictive than assuming strong convexity on the objective $f(\bx)$. 
	For example, assumption (e) is satisfied when $f(\bx)$ is convex and the convex function $\phi(\bx)$ defining the Bregman divergence is $\mu$-strongly convex.
\end{remark}

\begin{remark}
	Assumption (f) is satisfied if 
	%$f(\bx)$ is convex and $g(\bx\mid\bx_n) = f(\bx) + %
	%$B_{\phi}(\bx\Vert\bx_n)$, where the convex function
	$\phi(\bx)$ is $L$-smooth.
	Assumption (f) can be replaced by
	\begin{itemize}
		\item[(f\;\textprime)] $g(\bx \mid \by)$ is continuous in $\by$ in $D$.
	\end{itemize}
	This is the condition implicitly imposed in the proof of Proposition 7.4.1 in \citet{lange2016mm}. (This assumption is not made perfectly clear in the statement of the proposition.)
	Assumption (f\;\textprime) is satisfied, when 
	%if $f(\bx)$ is convex and $g(\bx\mid\bx_n) = f(\bx) + B_{\phi}(\bx\Vert\bx_n)$, where 
	$\phi(\bx)$ is a Bregman-Legendre function \citep{Byrne08Sequential,Byrne14lecture}.
\end{remark}

\subsection{Examples}

\subsubsection{Proximal gradient method}
The proximal gradient method minimizes $f(\bx) = f_0(\bx) + h(\bx)$ over $C=\mathcal{X}$, where both $f_0(\bx)$ and $h(\bx)$ are convex, proper, and lower semicontinuous. It is further assumed that $f_0(\bx)$ is $L$-smooth. The algorithm iteratively solves
\begin{equation}\label{eqn:proxgrad}
	\bx_{n+1} = \argmin_{\bx} \left\{ f_0(\bx_n) + \langle \nabla f_0(\bx_n), \bx - \bx_{n} \rangle + h(\bx) + \frac{1}{2\alpha}\|\bx - \bx_n\|^2 \right\}
\end{equation}
for a step size $0<\alpha < L^{-1}$.
To see that the proximal gradient algorithm is an instance of Bregman majorization, set $\phi(\bx) = \frac{1}{2\alpha}\|\bx\|^2 - f_0(\bx)$. Then 
\begin{equation}\label{eqn:bregma-proxgrad}
\begin{split}
	f(\bx) + B_{\phi}(\bx\Vert\bx_n) &= f_0(\bx) + h(\bx) + \frac{1}{2\alpha}\|\bx\|^2 - f_0(\bx) - \frac{1}{2\alpha}\|\bx_n\|^2 + f_0(\bx_n) \\
	& \quad - \langle \frac{1}{\alpha}\bx_n - \nabla f_0(\bx_n), \bx - \bx_n \rangle \\
%	&= f_0(\bx_n) + \langle \nabla f_0(\bx_n), \bx - \bx_n \rangle + h(x)
%		+ \frac{1}{2\alpha}\|\bx\|^2 + \frac{1}{2\alpha}\|\bx_n\|^2 - \frac{1}{\alpha}\langle \bx, \bx_n \rangle
%		\\
	&=  f_0(\bx_n) + \langle \nabla f_0(\bx_n), \bx - \bx_n \rangle + h(x)
		+ \frac{1}{2\alpha}\|\bx - \bx_n\|^2
\end{split}
\end{equation}
as desired.
It remains to verify that $f(\bx)$ and $\phi(\bx)$ satisfy conditions (a) through (f) of Propositions \ref{prop:objectiveconvergence} and \ref{prop:iterateconvergence}. Conditions (a) and (c) are assumed; (b) and (d) are true. Condition (e) is satisfied
since $\alpha \in (0, 1/L)$. The following fact is well-known:
\begin{lemma}
	A differentiable convex function $f(\bx)$ is $L$-smooth $\nabla f(\bx)$ if and only if $\frac{L}{2}\|\bx\|^2 - f(\bx)$ is convex.
\end{lemma}
\noindent Then, since $\phi(\bx) = \dfrac{1}{2}\left(\dfrac{1}{\alpha} - L\right)\|\bx\|^2 + \dfrac{L}{2}\|\bx\|^2 - f(\bx)$ and $\frac{1}{\alpha} > L$, $\phi$ is $(\frac{1}{\alpha} - L)$-strongly convex.

To check condition (f), we invoke the Baillon-Haddad theorem:
\begin{lemma}\label{lemma:Baillon-Haddad}
	If function $f(\bx)$ is convex, differentiable, and is $L$-smooth, then
$$
	\langle \nabla f(\bx) - \nabla f(\by), \bx - \by \rangle \ge \frac{1}{L}\|\nabla f(\bx) - \nabla f(\by) \|^2
	.
$$
\end{lemma}
\noindent Note $\nabla\phi(\bx) = \frac{1}{\alpha}\bx - \nabla f_0(\bx)$. Then,
\begin{align*}
	\|\nabla\phi(\bx) - \nabla\phi(\by)\|^2 &= \|\alpha^{-1}(\bx - \by) - [\nabla f_0(\bx) - \nabla f_0(\by)] \|^2 \\
	&= \frac{1}{\alpha^2}\|\bx - \by\|^2 + \|\nabla f_0(\bx) - \nabla f_0(\by)\|^2 
		- \frac{2}{\alpha}\langle \bx - \by, \nabla f_0(\bx) - \nabla f_0(\by) \rangle \\
	&\le \frac{1}{\alpha^2}\|\bx - \by\|^2 + \|\nabla f_0(\bx) - \nabla f_0(\by)\|^2 
		- \frac{2}{\alpha L}\| \nabla f_0(\bx) - \nabla f_0(\by) \|^2 \\
	&\le \frac{1}{\alpha^2}\|\bx - \by\|^2
	.
\end{align*}
The first inequality is due to Lemma \ref{lemma:Baillon-Haddad}. The last inequality holds since $\alpha \in (0, 1/L)$ implies $1 - \frac{2}{\alpha L} \le 0$.
Therefore $\nabla\phi(\bx)$ is $1/\alpha$-Lipschitz continuous and condition (f) is satisfied.

We summarize the discussion above as follows:
\begin{proposition}\label{prop:proxgrad}
	Suppose $f_0(\bx)$ and $h(\bx)$ are convex, proper, and lower semicontinuous. If $f_0(\bx)$ is $L$-smooth, then for $0 < \alpha < 1/L$, the proximal gradient iteration \eqref{eqn:proxgrad} converges to a minimizer of $f(\bx) = f_0(\bx) + h(\bx)$ if it exists.
\end{proposition}

\begin{remark}
	Lemma \ref{lemma:Baillon-Haddad} suggests that $\nabla\phi$ is $1/\alpha$-Lipschitz continuous if $0 < \alpha < 2/L$; in other words,  the step size may be doubled. Indeed, employing monotone operator theory \citep{Bauschke11convex,Ryu16primer} it can be shown that iteration \eqref{eqn:proxgrad} converges for $1/L \le \alpha < 2/L$ as well. Even though the MM interpretation is lost for this range of step size, the descent property remains intact \citep{She09thresholding,Bayram15convergence}.
\end{remark}
\begin{remark}
	The assumption that $h(\bx)$ is convex can be relaxed: if $h(\bx)$ is $\rho$-weakly convex, which means $h(\bx) + \frac{\rho}{2}\|\bx\|^2$ is convex, and $f_0(\bx)$ is $\rho$-strongly convex as well as $L$-smooth (this implies $\rho \le L$), then the objective $f(\bx)$ remains convex. The inner optimization problem in iteration \eqref{eqn:proxgrad} is also strongly convex if $\rho\alpha < 1$ and $\bx_{n+1}$ is uniquely determined. The latter condition is guaranteed if $\alpha \in (0, 1/L)$, and the conclusion of Proposition \ref{prop:proxgrad} holds. In fact, by using monotone operator theory, a larger step size $\alpha \in (0, \frac{2}{L+\rho})$ is allowed \citep{Bayram15convergence}. Statistical applications include nonconvex sparsity-inducing penalties such as the MCP \citep{Zhang10MCP}.
\end{remark}

\subsubsection{Mirror descent method}
For the constrained problem \eqref{eqn:constrained} and the Euclidean norm $\|\cdot\|_2$ the proximal gradient method takes the form of \emph{projected gradient}
\begin{equation}\label{eqn:projgrad}
\begin{split}
	\bx_{n+1} &= \argmin_{\bx\in C} \left\{ f(\bx_n) + \langle \nabla f(\bx_n), \bx - \bx_{n} \rangle + \frac{1}{2\alpha}\|\bx - \bx_n\|_2^2 \right\}
	\\
	&= P_C\left(\bx_n - \alpha \nabla f(\bx_n)\right)
	.
\end{split}
\end{equation}
This method relies heavily on the Euclidean geometry of $\mathbb{R}^d$, not $C$: $\|\cdot\|_2 = \langle \cdot, \cdot \rangle$.
If the distance measure $\frac{1}{2}\|\mathbf{x} - \mathbf{y}\|_2^2$ is replaced by something else (say $d(\mathbf{x}, \mathbf{y})$) that better reflects the geometry of $C$, then update such as
\begin{equation}\label{eqn:mirrordescent}
    \bx_{n+1} = P_C^{d}\left( \arg\min_{\bx\in \mathbb{R}^d} \left\{ f(\bx_n) + \langle \nabla f(\bx_{n}), \bx - \bx_n \rangle + \frac{1}{\alpha}d(\bx, \bx_{n}) \right\} \right)
\end{equation}
    may converge faster. Here, 
$$
    P_C^d(\mathbf{y}) = \argmin_{\bx\in C} d(\mathbf{x}, \mathbf{y})
$$
is a new (non-Euclidean) projection operator that reflects the geometry of $C$.

To see that iteration \eqref{eqn:mirrordescent} is a Bregman majorization for an appropriately chosen $d(\cdot,\cdot)$,  let
$$
    d(\bx, \by) = B_{\psi}(\bx\Vert\by) = \psi(\bx) - \psi(\by) - \langle \nabla \psi(\by), \bx - \by \rangle  \ge \frac{1}{2}\|\bx - \by\|^2
$$
for a $1$-strongly convex (with respect to some norm $\|\cdot\|$) and continuously differentiable function $\psi$ in $C$,
and set $\phi(\bx) = \frac{1}{\alpha}\psi(\bx) - f(\bx)$.
Similarly to equation \eqref{eqn:bregma-proxgrad}, we have
\begin{align*}
    f(\bx) + B_{\phi}(\bx\Vert\bx_n) 
	&=  f(\bx_n) + \langle \nabla f(\bx_n), \bx - \bx_n \rangle \\
	& \quad	+ \frac{1}{\alpha}\left[\psi(\bx) - \psi(\bx_n) - \langle \nabla\psi(\bx_n), \bx - \bx_n \rangle \right]
	\\
	&=  f(\bx_n) + \langle \nabla f(\bx_n), \bx - \bx_n \rangle + \frac{1}{\alpha} d(\bx, \bx_n)
	.
\end{align*}
Let $\tilde{\bx}_{n+1}$ be the \emph{unconstrained} minimizer of $f(\bx) + B_{\phi}(\bx\Vert\bx_n)$ (which is unique since $d(\bx, \bx_n)$ is strongly convex in $\bx$). The associated optimality condition entails
\begin{equation}\label{eqn:bregman_optimality}
    \nabla \psi(\tilde{\bx}_{n+1}) = \nabla\psi(\bx_{n}) - \alpha \nabla f(\bx_{n})
\end{equation}
Then,
\begin{align*}
    \bx_{n+1} &= \argmin_{\bx\in C}d(\bx, \tilde{\bx}_{n+1}) \\
        &= \argmin_{\bx\in C}\left\{ \psi(\bx) - \psi(\tilde{\bx}_{n+1}) - \langle \nabla\psi(\tilde{\bx}_{n+1}), \bx - \tilde{\bx}_{n+1} \rangle \right\} \\
        &= \argmin_{\bx\in C}\left\{ \psi(\bx) - \langle \nabla\psi(\tilde{\bx}_{n+1}), \bx \rangle \right\} \\
        &= \argmin_{\bx\in C}\left\{ \psi(\bx) - \langle \nabla\psi(\bx_{n}) - \alpha \nabla f(\bx_{n}), \bx - \bx_{n} \rangle \right\} 
        \\
        %&= \argmin_{\bx\in C}\left\{ \alpha f(\bx) + \alpha\phi(\bx) - \langle \alpha\nabla\phi(\bx_{n}), \bx - \tilde{\bx}_{n+1} \rangle \right\} 
        %\\
        &= \argmin_{\bx\in C}\left\{ f(\bx) + \phi(\bx) - \langle \nabla\phi(\bx_{n}), \bx - \bx_{n} \rangle  - \phi(\bx_n)  \right\} 
        \\
        %&= \argmin_{\bx\in C}\left\{ f(\bx) + B_{\phi}(\bx\Vert\bx_n) + \langle \nabla \phi(\bx_n), \bx_n - \tilde{\bx}_{n+1} \rangle \right\} 
        %\\
        &= \argmin_{\bx\in C}\left\{ f(\bx) + B_{\phi}(\bx\Vert\bx_n) \right\} 
        ,
\end{align*}
as sought.
To establish iterate convergence via SUMMA, we see that just as the proximal gradient method,  $f(\bx)$ and $\phi(\bx)$ satisfy conditions (a) through (e) of Propositions \ref{prop:objectiveconvergence} and \ref{prop:iterateconvergence} if $f$ is $L$-smooth and $\alpha \in (0, 1/L)$. In particular,
$$
    \phi(\bx) = \frac{1}{\alpha}\psi(\bx) - f(\bx) \ge \frac{1}{2\alpha}\|\bx\|^2 - f(\bx)
$$
to check condition (e). Condition (f\textprime) is fulfilled since $B_{\phi}(\bx\Vert\by)=\phi(\bx) - \phi(\by) - \langle \nabla\phi(\by), \bx - \by\rangle$ is continuous in $\by$ by construction.

Computation of $\bx_{n+1}$ can be further analyzed.
It is well known that if $\psi$ is $\mu$-strongly convex, then $\psi^*$ is $1/\mu$-smooth,
where $\psi^*$ is the Fenchel conjugate function of $\psi$:
$$
    \psi^*(\by) = \sup_{\bx\in\dom{\psi}} \langle \bx, \by \rangle - \phi(\bx)
$$
%(taking $\psi(\bx) = \infty$ if $\bx \notin C$).
\citep{Bauschke11convex}.
Hence $\nabla\psi^*$ is well-defined. 
Furthermore, $\nabla\psi^*(\nabla\psi(\bx))=\bx$.
%Since $\psi$ is strongly convex,  
Therefore the unconstrained optimality condition \eqref{eqn:bregman_optimality} is equivalent to
$$
    \tilde{\bx}_{n+1} = \nabla \psi^* \left(\nabla\psi(\bx_{n}) - \alpha \nabla f(\bx_{n}) \right),
$$
and we decompose the update \eqref{eqn:mirrordescent} into three steps:
\begin{align*}
    \by_{n+1} &= \nabla\psi(\bx_{n}) - \alpha \nabla f(\bx_{n}) \quad \text{(gradient step)} \\
    \tilde{\bx}_{n+1} &= \nabla \psi^*(\by_{n+1}) \quad \text{(mirroring step)} \\
    \bx_{n+1} &= P_C^{d}(\tilde{\bx}_{n+1}). \quad \text{(projection step)}
\end{align*}
Hence  Bregman majorization with $\phi(\bx)=\frac{1}{\alpha}\psi(\bx)-f(\bx)$ coincides with the mirror descent method under $B_{\psi}$ \citep{juditsky2011}.
The first step performs the gradient descent step in the dual space $\mathcal{X}^*$ of $\mathcal{X}$, 
and the second step maps the dual vector back to the primal space by the inverse mapping $\nabla\psi^*=(\nabla\psi)^{-1}$.
The final step projects (in a non-Euclidean fashion) the mapped primal vector onto the constraint set $C$.

\begin{example}[Exponentiated gradient]
    As a concrete instance of mirror descent, consider optimization over probability simplex $C = \Delta^{d-1} = \{\mathbf{x}\in\mathcal{X}=\mathbb{R}^d: \sum_{i=1}^d x_i = 1,~ x_i \ge 0, ~i=1, \dotsc, d\}$.
    An appropriate Bregman divergence is the Kullback-Leibler divergence, i.e., we use negative entropy $\psi(\bx) = \sum_{i=1}^d x_i\log x_i - \sum_{i=1}^d x_i$.
    %Pinsker's inequality states 
    It is easy to check, using the Taylor expansion and the Cauchy-Schwarz inequality,
    that $\psi$ is 1-strongly convex with respect to the $\ell_1$ norm $\|\bx\|_1=\sum_{i=1}^d |x_i|$ within $C$. 
    Furthermore, we have
    $\nabla\psi(\bx) = (\log x_1, \dotsc, \log x_d)^T =: \log\bx$
    and
    $\nabla\psi^*(\by) = (\nabla\phi)^{-1}(\by) = (e^{y_1}, \dotsc, e^{y_d})^T =: \exp(\by)$.
    The mirror descent or Bregmen MM update is then
    \begin{align*}
    \by_{n+1} &= \log\bx_{n} - \alpha\nabla f(\bx_{n}) \\
    %\tilde{\bx}_{n+1} &= \exp\left(\log\bx_{n} - \alpha\nabla f(\bx_{n}) \right) 
    \tilde{\bx}_{n+1} &= \exp\left(\by_n \right) 
    = \bx_{n}\odot\exp\left(- \alpha\nabla f(\bx_{n}) \right)
    \\
    %\bx_{n+1} &= \exp\left(\log\bx_{n} - \alpha\nabla f(\bx_{n}) \right)/Z_t 
    %= \bx_{n}\odot\exp\left(- \alpha\nabla f(\bx_{n}) \right)/Z_t,
    \bx_{n+1} &= \tilde{\bx}_{n+1} / Z_t,
    \end{align*}
    where $\odot$ denotes an elementwise product, and
    $$
    Z_t = \sum_{i=1}^d {x}_{n,i}\exp\left(- \alpha\nabla f(\bx_{n})_i \right) %= \mathbf{1}^T\exp\left(\log\bx_{n} - \alpha\nabla f(\bx_{n}) \right).
    $$
    is the normalization constant.
    The last step is because
    \begin{align*}
    P_{C}^{d}(\by) 
    &= \argmin_{\bx\in\Delta^{d-1}} B_{\psi}(\bx\Vert\by) 
    \\
    &= \argmin_{x_i \ge 0, \sum_{i=1}^d x_i=1} \sum_{i=1}^d\left(x_i\log\frac{x_i}{y_i}  - x_i + y_i\right)
    \\
    &= \argmin_{x_i \ge 0, \sum_{i=1}^d x_i=1} \sum_{i=1}^d\left(x_i\log\frac{x_i}{y_i}\right).
    \end{align*}
    and the associated Lagrangian 
$$
    \mathcal{L}(\bx, \mu) = \sum_{i=1}^d \left(x_i\log\frac{x_i}{y_i}\right) + \mu\left(\sum_{i=1}^d x_i - 1\right)
$$
    yields
$$
    %\log\frac{x_i}{y_i} + 1 + \mu = 0 
    %\quad \iff \quad
    x_i = y_i\exp(-\mu-1) = c y_i, %~ c > 0,
    \quad i=1, \dotsc, d.
$$
for some $c>0$.
    Summing these over all $i$ yields $c = 1/(\sum_{i=1}^d y_i)$ to have
$$
    x_i = \frac{y_i}{\sum_{j=1}^d y_j}, \quad i=1,\dotsc, d.
$$
    This special case is called the \emph{exponentiated gradient method} \citep{helmbold1997,azoury2001}.
\end{example}

\bibliographystyle{apalike}
\bibliography{bib-HZ,bib-JW}

\end{document}